\documentclass [11pt]{article}
\usepackage{mathtools,amsfonts,amssymb,amsthm, bm, nccmath}
\usepackage{epsfig}
\usepackage{graphicx}
\allowdisplaybreaks
\usepackage{enumerate}
\usepackage{cite}
\begin{document}
	\thispagestyle{empty}
	\null\vspace{-1cm}
	\medskip
	\vspace{1.75cm}
	\centerline{\Large\textbf{{Generalizations of Enestr\"{o}m-Kakeya Type}}}
		\centerline{\Large\textbf{{Theorems for matrix Polynomials}}}
	~~~~~~~~~~~~~~~~~~~~~~~~~~~~~~~~~~~~~~~~~~~~~~~~~~~~~~~~~~~~~~~~~~~~~~~~~~~~~~~~~~~~~~~~~~~~~~~~~~~~~~~~~~~~~~~~~~~~~~~~~~~~~~~~~~~~~~~~~~~~~~~~~~~~~~~~~~~

	\centerline{\bf {Idrees Qasim}}
	\centerline {Department of Mathematics, National Institute of Technology, Srinagar, India-190006}
	\centerline {idreesf3@nitsri.ac.in}
	~~~~~~~~~~~~~~~~~~~~~~~~~~~~~~~~~~~~~~~~~~~~~~~~~~~~~~~~~~~~~~~~~~~~~~~~~~~~~~~~~~~~~~~~~~~~~~~~~~~~~~~~~~~~~~~~~~~~~~~~~~~~~~~~~~~~~~~~~~~~~~~~~~~~~~~~~
	\vskip0.1in
	\noindent \textbf{Abstract}: In this paper, we establish bounds for the eigenvalues of matrix polynomials. Specifically, we find different generalizations of the Eneström-Kakeya Theorem for matrix polynomials.

	~~~~~~~~~~~~~~~~~~~~~~~~~~~~~~~~~~~~~~~~~~~~~~~~~~~~~~~~~~~~~~~~~~~~~~~~~~~~~~~~~~~~~~~~~~~~~~~~~~~~~~~~~~~~~~~~~~~~~~~~~~~~~~~~~~~~~~~~~~~~~~~~~~~~~~~~~~
	
	\noindent {{\bf Keywords:} Complex Polynomials, Matrix polynomial, Zeros, eigenvalues, Eneström-Kakeya Theorem.
	\vspace{0.15in}
	~~~~~~~~~~~~~~~~~~~~~~~~~~~~~~~~~~~~~~~~~~~~~~~~~~~~~~~~~~~~~~~~~~~~~~~~~~~~~~~~~~~~~~~~~~~~~~~~~~~~~~~~~~~~~~~~~~~~~~~~~~~~~~~~~~~~~~~~~~~~~~~~~~~~~~~~~
	
	\vspace{0.05in}
	~~~~~~~~~~~~~~~~~~~~~~~~~~~~~~~~~~~~~~~~~~~~~~~~~~~~~~~~~~~~~~~~~~~~~~~~~~~~~~~~~~~~~~~~~~~~~~~~~~~~~~~~~~~~~~~~~~~~~~~~~~~~~~~~~~~~~~~~~~~~~~~~~~~~~~~~~
	\section{Introduction} 
	Let $p(z):=\sum\limits_{j=0}^{m}a_jz^j$ represent a polynomial of degree $m$ with complex coefficients. The pursuit of finding zeros for such polynomials holds significance not only within mathematics but also across diverse fields like science, engineering, and technology. However, according to the Abel-Ruffini Theorem (refer to \cite{AR}), there exists no universal method for determining the zeros of a polynomial of degree 5 or higher using only arithmetic operations and radicals. Consequently, it becomes imperative to identify the regions where these zeros may lie. The problems involving polynomials in general, and location of their zeros in particular, have a long history and its application is in many areas of applied mathematics such as control theory, signal processing and electrical networks, coding theory and engineering among others \cite{16}.   In this context, Cauchy's classical theorem in 1829 (refer to \cite{CAU}) offers a valuable tool for analyzing the distribution of polynomial zeros, stating:\\
	
	\textbf{Theorem A.} All the zeros of the polynomial $p(z):=\sum\limits_{j=0}^{m}a_jz^j$, $a_m\neq 0$, with complex coefficients lie in
	$$\{z;|z|<s\}\subset \{z;|z|<1+M\},$$
	where $M:=\max\limits_{1\le j\le m-1}\left|\dfrac{a_j}{a_m}\right|$ and $s$ is the unique positive solution of
	$$|a_m|z^m-|a_{m-1}|z^{m-1}-\cdots-|a_1|x-|a_0|=0.$$
	
	Since then, numerous mathematicians have endeavored to establish improved bounds; for instance, see \cite{M,AM,MM,RS}.
	
	One of the elegant result in the theory of distribution of polynomials is the following theorem known as Enestr\"{o}m-Kakeya Theorem.\\
	
	\textbf{Theorem B.} If $p(z):=\sum\limits_{j=0}^{m}a_jz^j$ is a polynomial of degree $m$ with real coefficients satisfying $a_m\ge a_{m-1}\ge \dots\ge a_1\ge a_0>0,$ then all the zeros of $p$ lie in $|z|\le 1$.\\
	
	The theorem was proved by Enestr\"{o}m \cite{55}, independently by Kakeya \cite{85} and Hurwitz \cite{79}.\\
	
	There have been several generalizations and refinements of this result (for example, see \cite{AM, MM}).\\
	
	\noindent Consider $\mathbb{M}_{n,n}$ as the set comprising all $n\times n$ matrices with entries from $\mathbb{C}$. Herein, we define a matrix polynomial, denoted as $P:\mathbb{C}\rightarrow \mathbb{M}_{(n, n)}$, as a function expressed by
	$$P(z):=\sum_{j=0}^{m}A_jz^j,~~~A_j\in  \mathbb{M}_{(n, n)}.$$
	If $A_m\neq 0$, then $P(z)$ is termed a matrix polynomial of degree $m$. An eigenvalue of $P(z)$, denoted by $\lambda$, is a value for which there exists a non-zero vector $u\in\mathbb{C}^n$ such that $P(\lambda)u=0$. In such cases, $u$ is called an eigenvector of $P(z)$. The task of finding a number $\lambda\in \mathbb{C}$ and a non-zero vector $u\in \mathbb{C}^n$ satisfying $P(\lambda)u=0$ is referred to as a Polynomial Eigenvalue Problem ($PEP$). When $m=1$, this reduces to a Generalized Eigenvalue Problem ($GEP$) represented as $Au=\lambda Bu$. Moreover, if $B$ equals the identity matrix $I$, it becomes the standard eigenvalue problem $Au=\lambda u$.\\
	
	The computation of eigenvalues for a matrix polynomial poses a significant challenge. Iterative methods are commonly employed to tackle this challenge (for references, see \cite{11}). When computing pseudospectra of matrix polynomials, it becomes crucial to delineate a region within the complex plane encompassing the eigenvalues of interest. Bounds can aid in delineating such a region (for references, see \cite{7}). Moreover, bounding the eigenvalues can expedite numerical computations by aiding in selecting appropriate numerical methods and reducing the number of iterations required in various numerical algorithms.
	
	Note that if $A_0$ is singular, then 0 is an eigenvalue of $P(z)$, and if $A_m$ is singular, then 0 is an eigenvalue of the matrix polynomial $z^mP(1/z)$. Therefore, to locate the eigenvalues of these matrix polynomials, we always assume that $A_0$ and $A_m$ are non-singular.
	
	For a matrix $A\in \mathbb{M}_{(n,n)}$, $A\ge 0$ and $A>0$ is written if $A$ is positive semi-definite or positive definite respectively. Also for $B\in \mathbb{M}_{(n,n)}$, $A\ge B$ and $A>B$ means $A-B\ge 0$ and $A>B$ respectively. Moreover, $A^*$ denotes the transpose of $A$ and $tr(A)$ denotes trace of $A$. A vector $u\in \mathbb{C}^n$ is a unit vector if $||u||:=\sqrt{u^*u}=1,$ where $u^*$ denotes conjugate transpose of $u$.\\
	
	We note that any matrix-valued function $F(z)$ analytic in $|z|\le t$ can be expressed as a power series $F(z)=\sum\limits_{j=0}^{\infty}A_jz^j$, $A_j\in \mathbb{M}_{(n,n)},~~|z|\le t$ (for ref. see \cite{13}.).\\
	
	\noindent \textbf{Notation}. Throughout this paper, $||\cdot||$ denotes a subordinate matrix norm.\\

	The following extensions of Theorem A to matrix polynomial were obtained in \cite{CAM,7, 2}.
	
	\noindent \textbf{Theorem C.} Let $P(z):=\sum\limits_{j=0}^{m}A_jz^j$, $\det(A_m)\neq 0$ be a matrix polynomial. Then the eigenvalues of $P(z)$ lie in $|z|\le \rho$, where $\rho$ is a unique positive root of the equation
	
	$$||A_m^{-1}||^{-1}z^m-||A_{m-1}||z^{m-1}-\cdots-||A_1||-||A_0||=0.$$
	
	\noindent \textbf{Theorem D.} Let $P(z):=\sum\limits_{j=0}^{m}A_jz^j$, $\det(A_m)\neq 0$ be a matrix polynomial. Then the eigenvalues of $P(z)$ lie in $|z|<1+M$, where $M=||A_m^{-1}||\max\limits_{0\le j\le m-1}||A_j||.$
	
	Dirr and Wimmer \cite{DW} proved the following analogue of Theorem B concerning the bounds on the eigenvalues of matrix polynomials. 
	
	\noindent \textbf{Theorem E.} Let $P(z):=\sum\limits_{j=0}^{m}A_jz^j,$ be a matrix polynomial of degree $m$ such that
	$$A_m\ge A_{m-1}\ge \cdots \ge A_0\ge 0,~~A_m\neq 0,$$
	then all the eigenvalues of $P(z)$ lie in $|z|\le 1$.\\
	
	In this paper we find some results concerning bounds for the eigenvalues of matrix polynomials based on the norm of their coefficient matrices which gives generalizations of Theorem E with weakened condition.\\ 
	\section{Lemmas}
	For the proof of the theorems we need the following lemmas. The first Lemma is due to papescu et al. \cite{PJ}.\\
	\textbf{Lemma 1.} If $\mathcal{P}_k$ stands for the $k^{th}$ pell number defined by $\mathcal{P}_{0}=0,$ $\mathcal{P}_1=1$ and for $m\ge 2$, $\mathcal{P}_m=2\mathcal{P}_{m-1}+\mathcal{P}_{m-2}$, then
	$$\sum_{k=0}^{m}\binom{2m}{m+k}\mathcal{P}^2_k=2^{3(m-1)}.$$
	
	Next Lemma is due to Monga and Shah \cite{MS}.\\
	\textbf{Lemma 2.} Let $A,B, \in \mathbb{M}_{(n, n)}$ be such that $||A||_F\ge ||B||_F$ and $\angle (A,B)=\theta\le 2\alpha \le \pi$, then
	$$||A-B||_F\le (||A||_F-||B||_F)\cos \alpha+(||A||_F+||B||_F)\sin \alpha,$$
	where $||\cdot||_F$ is the Frobenius norm induced by the Frobenius inner product $\langle A \;,B \; \rangle=tr(B^*A)$. 
	
	We also need the following lemmas (for ref. see \cite{HJ}).\\
	\textbf{Lemma 3} Let $A\in \mathbb{M}_{(n, n)}$, then
	
	$$r(A)\le ||A||\le ||A||_F,$$
	where $r(A)=\max\{|u^*Au|; ||u||=1\}$ is called as the numerical radius of $A$.\\
	
	\noindent \textbf{Lemma 4.} Let $A\in \mathbb{M}_{n, n}$, be a Hermitian matrix, then 
	$$\lambda_{\min}(A)=\min_{||u||=1}\{u^*Au\}\le \max_{||u||=1}\{u^*Au\}=\lambda_{\max}(A).$$

	\section{Main Results}
	The following theorem gives generalization of Theorem B for matrix polynomials.\\
\textbf{Theorem 1.} Let $P(z):=\sum_{j=0}^{m}A_jz^j$ is a matrix polynomial of degree $m$, where $A_m$ is invertible. For $t>0$, let 
$$||A_m^{-1}||^{-1}\ge t||A_{m-1}||\ge t^2||A_{m-2}||\ge\cdots\ge t^m||A_0||,~~A_0>0,$$
then $P(z)$ has all its eigenvalues in $|z|\le k_1/t$, where $k_1$ is the greatest positive root of the trinomial equation $K^{m+1}-2K^m+1=0$.
\begin{proof}
Define
$$g(z)=A_{m-1}z^{m-1}+A_{m-2}z^{m-2}+\cdots+A_1z+A_0,$$
and let $u$ be a unit vector, then we have for $|z|=R(>1/t),$
\begin{align*}
	||g(z)u||&\le ||A_{m-1}||R^{m-1}+||A_{m-2}||R^{m-2}+\cdots +||A_0||\\
	&=||A_{m-1}||R^{m-1}\bigg\{1+\frac{||A_{m-2}||}{||A_{m-1}||}\frac{1}{R}+\frac{||A_{m-3}||}{||A_{m-1}||}\frac{1}{R^2}+\cdots+\frac{||A_{1}||}{||A_{m-1}||}\frac{1}{R^{m-2}}\\
	&~~~~~~~~~~~~~~~~~~~~~~~~~~~~~~~~~+\frac{||A_{0}||}{||A_{m-1}||}\frac{1}{R^{m-1}}\bigg\}\\
	&\le||A_{m-1}||R^{m-1}\left\{1+\frac{1}{tR}+\frac{1}{(tR)^2}+\cdots+\frac{1}{(tR)^{m-2}}+\frac{1}{(tR)^{m-1}}\right\}\\
	&=||A_{m-1}||\left\{\frac{(tR)^m-1}{t^{m-1}(tR-1)}\right\}.
\end{align*}	
Thus, we have for $|z|=R$,
\begin{align*}
	||P(z)u||&=||A_mz^mu+g(z)u||\\
	&\ge ||A_mz^mu||-||g(z)u||\\
	&\ge ||A_m^{-1}||^{-1}R^m-||g(z)u||\\
	&\ge||A_m^{-1}||^{-1}R^m-||A_{m-1}||\left\{\frac{(tR)^m-1}{t^{m-1}(tR-1)}\right\}\\
	&>0,
\end{align*}
if
$$\frac{||A_m^{-1}||^{-1}}{t||A_{m-1}||}>\frac{(tR)^m-1}{(tR)^{m}(tR-1)}.$$
That is, if
$$1>\frac{(tR)^m-1}{(tR)^{m}(tR-1)}.$$
That is, if
\begin{equation*}
	(tR)^{m+1}-2(tR)^m+1>0.
\end{equation*}
This implies that all the eigenvalues of $P(z)$ lie in $|z|\le k_1/t$, where $k_1$ is the greatest positive root of the trinomial equation
$$K^{m+1}-2K^m+1=0.$$
This proves the theorem.
\end{proof}
\noindent	\textbf{Theorem 2.} Let $P(z):=\sum\limits_{j=0}^{m}A_jz^j$ is a matrix polynomial of degree $m$, such that for $k\ge 1$, $k||A_m||_F\ge ||A_{m-1}||_F\ge \cdots\ge ||A_1||_F\ge ||A_0||_F$ and $\angle(A_j,C)\le \alpha \le \dfrac{\pi}{2},~j=0,1,2,\dots,m$ for some non-zero matrix $C\in \mathbb{M}_{n, n}$. Then all the eigenvalues of $P(z)$ lie in
\begin{equation}
	|z+k-1|\le \frac{1}{||A_m^{-1}||_F^{-1}}\left\{(k||A_m||_F-||A_0||_F)(\sin \alpha+\cos \alpha)+||A_0||_F+2\sin \alpha \sum_{j=0}^{m-1}||A_j||_F\right\}.\label{eq:01}
\end{equation}
\begin{proof}
	Consider the polynomial
	\begin{align*}
		F(z)&=(1-z)P(z)\\
		&=-A_mz^{m+1}+\sum_{j=0}^{m}(A_j-A_{j-1})z^j,~~~~A_{-1}=0.
	\end{align*}
Let $u$ be a unit vector, then we have
\begin{align*}
	||F(z)u||&\ge||(A_mz^{m+1}+kA_mz^m-A_m)u||-\left|\left|(kA_m-A_{m-1})z^m+\sum_{j=0}^{m-1}(A_j-A_{j-1})z^j\right|\right|\\
	&\ge |z|^m\left[||A_m(z+k-1)u||-\left\{||kA_m-A_{m-1}||+\sum_{j=0}^{m-1}\frac{||A_j-A_{j-1}||}{|z|^{m-j}}\right\}\right]\\
	&\ge |z|^m\left[||A_m^{-1}||^{-1}|z+k-1|-\left\{||kA_m-A_{m-1}||+\sum_{j=0}^{m-1}\frac{||A_j-A_{j-1}||}{|z|^{m-j}}\right\}\right].
\end{align*}
Using Lemma 3, it follows that
$$||F(z)u||\ge\left[||A_m^{-1}||_F^{-1}|z+k-1|-\left\{||kA_m-A_{m-1}||_F+\sum_{j=0}^{m-1}\frac{||A_j-A_{j-1}||_F}{|z|^{m-j}}\right\}\right].$$
Let $|z|>1$, so that $\dfrac{1}{|z|^{m-j}}<1$, $0\le j<m$ and using Lemma 2, we get
\begin{fleqn}
\begin{align*}
	||F(z)u||&\ge 
	|z|^m\bigg[||A_m^{-1}||_F^{-1}|z+k-1|-\bigg\{(k||A_m||_F-||A_0||_F)(\cos \alpha+\sin \alpha)+||A_0||_F\\
&~~~~~~~~~~~~~~~~~~~~~~~~~~~~~~~~~~~~~~~~~~~~~~~	+2\sin \alpha\sum_{j=0}^{m-1}||A_j||_F\bigg\}\bigg]\\
	&>0,
\end{align*}
\end{fleqn}
if
$$|z+k-1|> \frac{1}{||A_m^{-1}||_F^{-1}}\left\{(k||A_m||_F-||A_0||_F)(\sin \alpha+\cos \alpha)+||A_0||_F+2\sin \alpha \sum_{j=0}^{m-1}||A_j||_F\right\}.$$
This shows that the eigenvalues of $F$ having modulus greater than one lie in \eqref{eq:01}. But the eigenvalues of $F$ of modulus not greater than one already lie in \eqref{eq:01}. Since the eigenvalues of $P$ are also eigenvalues of $F$, the result follows.
\end{proof}

\noindent \textbf{Remark.} A result of Shah and Liman \cite{SL} follows when we take $n=1$.\\

When we apply Theorem 2 to $P(tz)$ we obtain the following.

\noindent \textbf{Corollary 1.} Let $P(z):=\sum_{j=0}^{n}A_jz^j$ such that for some $k\ge 1$ and $t>0$,
$$kt^n||A_m||_F\ge t||A_{m-1}||_F\ge \cdots\ge t||A_1||_F\ge ||A_0||_F$$
and $\angle(A_j,C)\le \alpha \le \dfrac{\pi}{2},~j=0,1,2,\dots,m$ for some non-zero matrix $C\in \mathbb{M}_{n, n}$. Then all the eigenvalues of $P(z)$ lie in
$$|z+kt-t|\le \frac{t}{||A_m^{-1}||_F^{-1}}\bigg[\left(k||A_m||_F-\frac{||A_0||_F}{t^m}\right)(\sin \alpha +\cos \alpha)+\frac{||A_0||_F}{t^m}~~~~~~~~~~~$$
$$~~~~~~~~~~~~~~~~~~~~~~~~~~~~~~~~~~~~~~~~~~~~~~~~~~~~~~~~~~~~~~~~~~+2\sin \alpha \sum_{j=0}^{m-1}||A_j||t^{j-m}\bigg].$$

Next, we use Pell number to find a ring shaped region containing all the eigenvalues of a matrix polynomial.\\
\noindent \textbf{Theorem 3.} Let $P(z):=\sum_{k=0}^{m}A_kz^k$, be a matrix polynomial. Then all its eigenvalues lie in the ring shaped region $D:=\{z\in \mathbb{C}: r_1\le |z|\le r_2\},$ where
\begin{equation}
	r_1=\min_{1\le k\le m}\left\{2^{3(1-m)}\binom{2m}{m+k}\mathcal{P}^2_k\frac{1}{||A_0^{-1}||||A_k||}\right\}^{1/k}\label{eq:031}
\end{equation}
and 
\begin{equation}
	r_2=\max_{1\le k\le m}\left\{2^{3(m-1)}\binom{2m}{m+k}^{-1}\mathcal{P}^{-2}_k||A_{m-k}||||A_m^{-1}||\right\}^{1/k}.\label{eq:032}
\end{equation}
\begin{proof}
	From equation \eqref{eq:031}, it follows that for $1\le k\le m$,
	\begin{align}
		& r_1^k\le 2^{3(1-m)}\binom{2m}{m+k}\mathcal{P}^2_k\frac{1}{||A_0^{-1}||||A_k||}\nonumber \\
		\Rightarrow~&||A_0^{-1}||||A_k||r_1^k\le 2^{3(1-m)}\binom{2m}{m+k}\mathcal{P}^2_k.\label{eq:033}
	\end{align}
	Let $u\in \mathbb{C}^n$ be the unit eigenvector. Assume that $|z|<r_1$, then
	\begin{align*}
		||P(z)u||&\ge ||A_0u||-||\sum_{k=1}^{m}A_kz^k u||\\
		&\ge ||A_0^{-1}||^{-1}-\sum_{k=1}^{m}||A_k|||z^k|\\
		&>||A_0^{-1}||^{-1}-\sum_{k=1}^{m}||A_k||r_1^k\\
		&=||A_0^{-1}||^{-1}\left(1-||A_0^{-1}||\sum_{k=1}^{m}||A_k||r_1^k\right)\\	
		&\ge ||A_0^{-1}||^{-1}\left(1-\sum_{k=1}^{m}2^{3(1-m)}\binom{2m}{m+k}\mathcal{P}^2_k\right)=0
	\end{align*}
on account of Lemma 1 and consequently $P(z)$ does not have any eigenvalue in $|z|<r_1$.\\
It is known (Theorem C) that all the eigenvalues of $P(z)$ have modulus less or equal than the unique positive root of the equation
$$G(z)=z^m||A_m^{-1}||^{-1}-z^{m-1}||A_{m-1}||-\cdots-z||A_1||-||A_0||=0.$$
Hence the second part of our result will be proved if we show that $G(r_2)\ge 0$.
From equation \eqref{eq:032}, it follows for $1\le k\le m$,
$$||A_m^{-1}||||A_{m-k}||\le 2^{3(1-m)}\binom{2m}{m+k}\mathcal{P}^2_kr^k_2.$$
Then
\begin{align*}
	G(r_2)&=r_2^m||A_m^{-1}||^{-1}-r^{m-1}_2||A_{m-1}||-\cdots-r_2||A_1||-||A_0||\\
	&=||A_m^{-1}||^{-1}\left(r^m_2-||A_m^{-1}||\sum_{k=1}^{m}||A_{m-k}||r_2^{m-k}\right)\\
	&\ge ||A_m^{-1}||^{-1}\left[r^m_2-\left(\sum_{k=1}^{m}2^{3(1-m)}\binom{2m}{m+k}\mathcal{P}^2_kr^k_2\right)r^{m-k}_2\right]\\
	&=||A_m^{-1}||^{-1}r^m_2\left[1-\left(\sum_{k=1}^{m}2^{3(1-m)}\binom{2m}{m+k}\mathcal{P}^2_k\right)\right]\\
	&=0,
\end{align*}
and the proof is complete.
\end{proof}
Next, we prove the following generalization of a result due to Monga and Shah \cite{MS}.\\

\noindent \textbf{Theorem 4.} If $f(z):=\sum_{j=0}^{\infty}A_jz^j$ is analytic in $|z|<t$, such that for some $k\ge 1$,
$$kA_0\ge tA_1\ge t^2A_2\ge \cdots,$$
then $f(z)$ does not vanish in
$$\left|z-\left(\frac{k-1}{2k-1}\right)t\right|\le \frac{kt}{2k-1}.$$
\begin{proof}
	Since $f(z):=\sum_{j=0}^{\infty}A_jz^j$ is analytic in $|z|\le t$, therefore $\lim\limits_{j\rightarrow \infty}A_jz^j=0$. Let $u$ be a unit vector. Define $F_u(z)=u^*(z-t)f(z)u$, then $F_u(z)$ is a complex function analytic in $|z|\le t$. Also,
	$$u^*kA_0u\ge u^*tA_1u\ge u^*t^2A_2u\ge \cdots.$$
	Now,
	\begin{align}
		F_u(z)&=-u^*tA_0u+u^*(A_0-tA_1)zu+u^*(A_1-tA_2)z^2u+\dots\nonumber\\
		&=-u^*tA_0u+u^*A_0zu-u^*kA_0zu+u^*(kA_0-tA_1)zu+u^*z\sum_{j=2}^{\infty}(A_{j-1}-tA_j)z^{j-1}u\nonumber\\
		&=-u^*tA_0u+u^*A_0zu-u^*kA_0zu+G_u(z),\label{eq:BB01}
	\end{align} 
where
$$G_u(z)=u^*(kA_0-tA_1)zu+u^*z\sum_{j=2}^{\infty}(A_{j-1}-tA_j)z^{j-1}u.$$
Clearly $G_u(z)$ is analytic, $G_u(0)=0$ and for $|z|=t$
\begin{align*}
	|G_u(z)|& \le u^*(kA_0-tA_1)tu+t\sum_{j=2}^{\infty}u^*(A_{j-1}-tA_j)ut^{j-1}\\
	&=u^*kA_0tu.
\end{align*}
Therefore, by Schwarz's Lemma
\begin{equation}
	|G_u(z)|\le u^*kA_0u|z|~~~~\text{for}~~|z|<t.
\end{equation}
Thus, it follows from \eqref{eq:BB01}
\begin{align*}
	|F_u(z)|&\ge |u^*(tA_0-A_0z+kA_0z)u|-|G_u(z)|\\
	&\ge |u^*A_0u|\left\{|(k-1)z+t-k|z|\right\}\\
	&>0, 
\end{align*}
if
$$k|z|<|(k-1)z+t|.$$
It is easy to verify that the region defined by
$$\left\{z:k|z|<|(k-1)z+t|\right\}$$
is precisely the disk
\begin{equation}
	\left\{z:\left||z-\left(\frac{k-1}{2k-1}\right)t\right|<\frac{kt}{2k-1}\right\}.\label{eq:MBH1}
\end{equation}
It follows that $F_u(z)$ and hence $f(z)$ does not vanish in the disk defined by (\ref{eq:MBH1}). This proves the theorem.
\end{proof}

\section{Declaration}
\begin{itemize}
	\item Availability of data and materials: Not applicable
	\item Competing interests: The author declare that we have no competing interest.
	\item Funding: The author does not receive any funding for this research
	\item Authors' contributions - provide individual author contribution: Not applicable
	\item Acknowledgements: 
	\item Authors' information: Department of Mathematics, National Institute of Technology, Srinagar, India, 190006. Email: idreesf3@nitsri.ac.in
\end{itemize}
	
\end{document}